\documentclass{amsart}

\newtheorem{thm}{Theorem}[section]
\newtheorem{ass}[thm]{Assumption}

\newtheorem{lem}[thm]{Lemma}
\newtheorem{prop}[thm]{Proposition}

\theoremstyle{definition}

\newtheorem{eg}[thm]{Example}

\theoremstyle{remark}

\newtheorem{remk}[thm]{Remark}

\newcommand{\Card}{{\rm Card}}

\usepackage{graphicx}
\usepackage{psfrag}

\begin{document}

\begin{abstract}We consider the spectral behavior of the Laplace-Beltrami operator
      associated to a class of singular perturbations of a Riemannian
      metric on a complete manifold.  The class of perturbations
      generalizes the well-known `opening node' perturbation of
      Teichm\"uller theory.  In particular, we recover results
      of Ji and Zworski \cite{JiZwr93} and Wolpert \cite{Wlp92}
      from our more general methods.
\end{abstract}

\dedicatory{For Ralph Phillips}         

\title[Tracking eigenvalues to the frontier]{
\mbox{Tracking eigenvalues to the 
      frontier of moduli space I:} 
      Convergence and spectral accumulation}

\author{Christopher M. Judge}       

\address{Indiana University, Bloomington, IN.}  

\email{cjudge@indiana.edu, 
          http://www.math.indiana.edu}

\thanks{This work was supported in part by 
        the National Science Foundation under
        DMS-9972425}

\subjclass{58G25, 35P20}

\keywords{Laplacian, degeneration, 
         hyperbolic surface, eigenvalue}

\date{\today}

\maketitle

\section{Introduction}

\begin{figure}
  \psfrag{t}{$t \in I$}
  \psfrag{M}{$N$}
  \psfrag{U}{$I \times M$}
  \centering
    \includegraphics{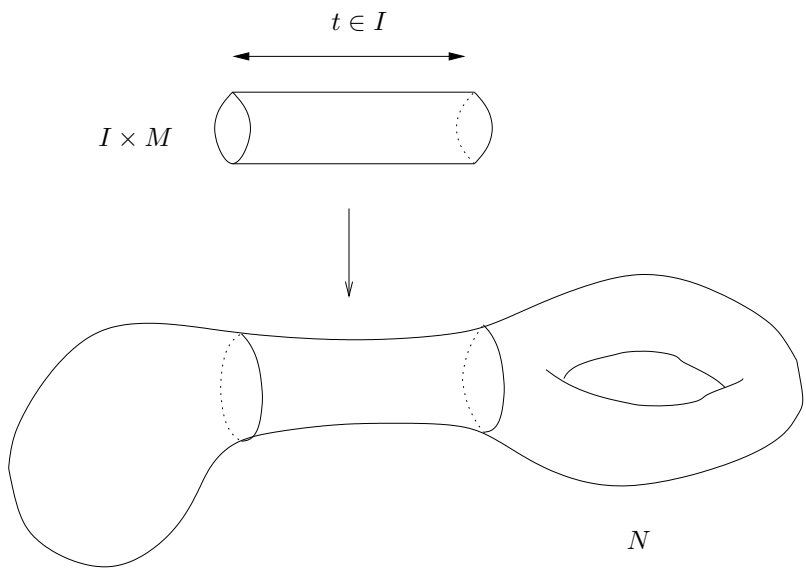}
  \caption{The set up}
  \label{Degeneration}
\end{figure}

Let $M$ and $N$ be compact 
differentiable manifolds of dimension $d$
and $d+1$ respectively. Given a compact
interval $I \subset {\mathbb R}$, 
we suppose that we have a (fixed) 
embedding $I \times M \subset N$. 
See Figure \ref{Degeneration}.
Let $N^0$ denote the complement of 
$\{0\} \times M$.

Let $h$ be a Riemannian metric on $M$,
and let $\rho$ be a positive function 
on ${\mathbb R}^2$ that is positively 
homogeneous of degree $1$ and smooth 
away from $(0,0)$. Given
$(a,b) \in {\mathbb R^2}$, we
consider continuous families, $g_{\epsilon}$,
of symmetric $(0,2)$-tensors on $N$
such that 
\begin{equation} \label{LocalMetric}
  g_{\epsilon} |_{I \times M}~ =~
 \rho(\epsilon,t)^{2a} \cdot dt^2~ 
 +~ \rho(\epsilon,t)^{2b} \cdot h
\end{equation}
and such that $g_{\epsilon}$ is positive 
definite on $N^0$. Note that 
for $\epsilon \neq 0$, the tensor 
$g_{\epsilon}$ is a Riemannian metric 
on all of $N$, but that the tensor
$g_0$ is singular on 
$\{0\} \times M$ if $(a,b) \neq (0,0)$.

\medskip

\begin{eg}[Hyperbolic Degeneration]
\label{DegenExample}
Let $\gamma$ be a simple closed curve on 
a compact oriented surface $N$ with $\chi(N) <0$. 
Let $g_{\epsilon}$ be a metric on $N$
of constant curvature $-1$ 
such that the unique geodesic 
homotopic to $\gamma$ has length 
$\epsilon< 2 \cosh^{-1}(2)$. 
By the collar lemma \cite{Bsr},
there exists an embedding 
$I \times \gamma \rightarrow N$ with
$I=[-1, 1]$ such that\footnote{
    The coordinates given here for 
    a collar made their first 
    appearance in \cite{JdgPhl97}.
 To obtain the more common Fermi
 coordinates, let 
  $t=\epsilon \sinh(\rho)$.}
\begin{equation} \label{Classical}
  g_{\epsilon}|_{I\times \gamma}~ =~ 
  \frac{dt^2}{\epsilon^2 +t^2} 
  +  (\epsilon^2 +t^2)~ d x^2
\end{equation}
where $x$ is the usual coordinate on
the circle 
${\mathbb R} / {\mathbb Z} \cong \gamma$.
Note that the Riemannian surface 
$((I \times \gamma)^0,g_0)$ is 
a union of hyperbolic cusps.
\end{eg}

\medskip

In both this paper and 
\cite{Jdg00} we study 
the small $\epsilon$ behavior 
of the spectrum of the 
Laplace-Beltrami operator 
$\Delta_{\epsilon}$ associated to a
metric families $g_{\epsilon}$  as 
described above in (\ref{LocalMetric}). 
{\em In this paper,  we restrict our
attention to $(a,b)$ satisfying
$a \leq -1$ and $b>0$}.

\medskip

\begin{thm}[Theorem 
\ref{EigenfunctionConvergence}] 
\label{Main1}
Let $\epsilon_j \rightarrow 0$. 
Any sequence $\psi_{j}$ of eigenfunctions
of $\Delta_{\epsilon_j}$ with uniformly 
bounded eigenvalues has a subsequence
that converges (up to rescaling) 
to an eigenfunction $\psi_* \neq 0$ 
of $\Delta_0$.
In particular, for each compact
subset of $A  \subset N^0$
the subsequence converges to $\psi_*$ in 
$H^1(A, dV_0)$.
\end{thm}

\medskip 

\noindent
In the special case of Example 
\ref{DegenExample}, the preceding 
theorem was obtained by Wolpert 
\cite{Wlp92} and Ji \cite{Ji93}.

As observed in \cite{Mlr} \S 8.1,
the manifold $(N^0, g_0)$ is 
Riemannian complete if and only if 
$a \leq -1$. In  
\S \ref{SectionManifoldCusps},
we prove that if $g_0$ is 
`marginally complete', $a=-1$, 
then the essential spectrum of 
$\Delta_0$ consists of the band 
$[(2^{-1}b d c)^2, \infty[$ where 
$c$ is determined by $\rho$.
See Proposition \ref{ContinuousSpectrum}.
In the `overcomplete' case, $a<-1$,
we show that the essential spectrum 
consists of the band $[0, \infty[$.

On the other hand, for  $\epsilon \neq 0$, 
the operator $\Delta_{\epsilon}$ 
has purely discrete spectrum. Hence one is
led to ask about the nature of the transition
from discrete to continuous spectrum as 
$\epsilon$ tends to zero. 
The following gives a precise 
quantitative description of this transition.

\medskip

\begin{thm}[Theorem \ref{CountTheorem}]
 \label{Main2}
Let $a=-1$, $b>0$, and let 
$c_{\pm}$ satisfy 
$\rho(0,t) = c_{\pm} t$ for $ \pm t>0$.
For $\Lambda>0$, let 
${\mathcal N}_{\Lambda}(\epsilon)$ 
be the number of eigenvalues of 
$\Delta_{\epsilon}$ that lie 
in $[0, \Lambda]$. Then 
\begin{equation*}  
%\label{MainEigenvalueCount1}
%
%
  {\mathcal N}_{\Lambda}(\epsilon) =
   \left( c_+
 \sqrt{\Lambda-
  \left(\frac{c_{+} \cdot b \cdot d}{2}
  \right)^2 } 
  +   c_{-} \sqrt{\Lambda- 
 \left( \frac{c_{-}\cdot b \cdot d}{2}
  \right)^2 }  \right) \cdot
   \frac{\log( \epsilon^{-1})}{\pi} 
   + O_{\Lambda}(1).
\end{equation*}
\end{thm}

\medskip

\noindent
L. Ji and M. Zworski \cite{JiZwr93}
obtained Theorem \ref{Main2} in the 
special case of  Example 
\ref{DegenExample} following the earlier 
work of  \cite{Ji93}  
and \cite{Wlp87}.

We now provide an outline 
of the present paper.
In \S \ref{SectionDegeneration} 
we establish notation.  
In \S \ref{SectionConstantMode} 
we separate the action of the 
Laplacian on functions that 
are constant on each fibre of 
$I \times M \rightarrow I$ from its 
action on those functions whose integral along
each fibre vanishes. Using a classical
transformation of Sturm-Liouville 
theory,
 we demonstrate a unitary equivalence 
of the former action with the action of
\begin{equation*} 
   \frac{\partial^2}{\partial s^2}~ 
   - \frac{bd}{2} \cdot (\alpha^{-1})^*
   \left( \rho^{-2a-2} \left(
   \rho \cdot \rho_{tt}''
   +\frac{bd-2a-2}{2}(\rho'_t)^2 \right)
   \right)
\end{equation*}
on $L^2(\alpha(I), ds)$ 
 where 
$\alpha$ is a diffeomorphism of intervals
determined by $\rho$ (Lemma \ref{Unitary}).
We also provide a convexity estimate on 
the fibrewise $L^2(M, dV_h)$-norm of  Laplace 
eigenfunctions whose integral along each
fibre is zero (Lemma \ref{ConvexLemma}).

In \S \ref{SectionManifoldCusps} we
obtain the basic spectral theory of 
a manifold with (generalized) cusps
by specializing the results
of \S \ref{SectionConstantMode}
to the case $\epsilon=0$.
For example, we show that the `cut-off' 
Laplacian $\Delta^{\perp}_0$ is compactly
resolved, and that the fibrewise 
$L^2(M, dV_h)$-norm of a
cusp form vanishes to arbitrary order 
at $t=0$.

In \S \ref{SectionEigenfunctionConvergence} 
we prove Theorem \ref{Main1} using the results
of earlier sections.
In \S \ref{SectionCutOffLaplacian} 
we show that the spectrum of
$\Delta^{\perp}_{\epsilon}$ varies 
continuously in $\epsilon$ 
including at $\epsilon=0$. 
We use this continuity in \S \ref{SectionCount}
together with Dirichlet-Neumann bracketing
and elementary Sturm-Liouville theory
to prove Theorem 
\ref{Main2}.

%%%%%%%%%%%%%%%%%%%%%%%%%%%%%%%%%%%%%%%%%%%%%%%%%%%%%%%%%%%%%%%%%%

\section{Preliminaries and notation}
\label{SectionDegeneration}

Let $g$ be a Riemannian metric 
on a manifold $N$, and let $dV_g$
and $\nabla_g$ denote respectively the 
associated density and gradient.\footnote{We 
will drop the subscript $g$ if it is clear 
from the context.} 
The Laplace-Beltrami operator $\Delta_g:L^2(N,dV_g) \rightarrow L^2(N,dV_g)$ 
is defined via the Friedrichs extension\footnote{See, for example,  
\cite{Kat} page 325.}  of the form 
$Q(\phi, \psi) = \int_{N} 
  g(\nabla \phi, \nabla \psi) dV_g$ 
with respect to the $L^2(N,dV_g)$-norm. 
In particular, by design
\begin{equation}   \label{BasicParts}
    \int_N g(\nabla_g \phi , \nabla_g \psi)~
  dV_g 
  = \int_N \phi \cdot \Delta_g \psi~dV_g
\end{equation}
for any $\phi$ and $\psi$ 
in the domain of $\Delta_g$.
The core used in this Friedrichs extension
is the space of smooth functions $f$ such 
that $Q(f,f) + \int_N f^2 < \infty$ satisfying
symmetric boundary conditions.

Consider a metric $g$ on $I \times M$
of the form given in (\ref{LocalMetric}), and
recall that $d$ denotes the dimension of $M$.
For any $f \in C^{\infty}_0(I \times M)$  
\begin{equation} \label{Laplacian}
 \Delta_g f = - L(f)
 + \rho^{-2b} \Delta_h f 
\end{equation}
where 
\begin{equation}
 \label{L}
  L(f) = \rho^{-a-bd}~
             \partial_t~  \rho^{-a+bd}~
       \partial_t~ f.
\end{equation}
and
\begin{equation}
 \label{Gradient}
 \nabla_g f = 
 \rho^{-2a} \cdot \partial_t f
 + \rho^{-2b} \cdot \nabla_h f .
\end{equation}
The volume form restricted to $I \times M$ is 
\begin{equation}
  \label{Volume}
   dV_g = \rho^{a+bd}~dt~dV_h.
\end{equation}

The function $t \rightarrow \rho(0,t)$ 
is positive and  positively homogeneous
of degree 1. Hence
there exist unique positive 
constants $c_+$ and $c_-$ such that 
\begin{equation} \label{Homo}
 \rho(0,t)= \left\{ 
\begin{array}{ll}
  c_+ t, &  t>0 \\
  c_- t,  &  t <0. 
\end{array}  \right.
\end{equation}
We will refer to $c_{\pm}$ as 
{\em constants of homogeneity}.

\medskip

%%%%%%%%%%%%%%%%%%%%%%%%%%%%%%%%%%%%%%%%%%%%%%%%%%%%%%%%%%%%%%

\section{The constant mode and its complement}

\label{SectionConstantMode}

\medskip

\begin{remk} 
 In this section, if $\epsilon=0$, 
 then replace $I$ with $I \setminus \{0\}$.
\end{remk}

\medskip

Viewing $I \times M$ as an $M$-fibre bundle 
over $I$, we use integration
along the fibres to analyse Laplace 
eigenfunctions on $I \times M$. 
For $f \in L^2(I \times M)$ define the 
{\em constant mode (zeroth Fourier coefficient)
of $f$} to be 
\begin{equation} 
   f_0(t) = \int_{\{t\}\times M} f(t,m)~ dV_h.
\end{equation}
One may regard $f_0$ as a function 
on $I \times M$ that is constant along
the fibres $\{t\} \times M$. From this
point of view $P_0(f)=f_0$ defines an
orthogonal projection onto a closed subspace
of $L^2(I \times M, dV_g)$.  
Let $P^{\perp}$ be the orthogonal projection
onto the complementary subspace 
$L^2_{\perp}(I \times M, dV_g)$.

These projections diagonalize the Laplacian:
\begin{equation} \label{Diagonal}
\Delta_g  =  P_0 \circ \Delta_g \circ P_0  + 
     P_{\perp} \circ \Delta_g \circ P_{\perp}.
\end{equation}
The operator $P_0~ \Delta_g~ P_0$ is 
unitarily equivalent to the operator 
$-L$ (densely) 
defined on $L^2( I, \rho^{a+bd} dt)$ 
via the Friedrichs extension. 
In particular, if $\Delta \psi= \lambda \psi$, 
then 
\begin{equation} \label{ConstantODE}
   -L~ \psi_0=  \lambda \cdot  \psi_0.
\end{equation}
The operator  $P_{\perp} \Delta_g P_{\perp}$ 
is unitarily equivalent to the restriction, 
$\Delta^{\perp}_{\epsilon}$, 
of $\Delta_{g_{\epsilon}}$ to 
$L^2_{\perp}(I \times M, dV_g)$.

To analyse solutions to (\ref{ConstantODE}),
we conjugate $L$ with a unitary operator
found in classical Sturm-Liouville theory 
(see, for example, \cite{CrnHlb} V \S 3.3).
For $t_{0} \in I$ and $\epsilon$ fixed, let
\begin{equation} \label{Alpha}
  \alpha(t)
  = \int_{t_0}^t \rho^{a}(\epsilon,u) du.
\end{equation}
and define 
$U: C^{\infty}(\alpha(I)) 
  \rightarrow C^{\infty}(I)$
by
\begin{equation} \label{UnitaryDefn}
 U(f) = 
   \rho^{-\frac{bd}{2}} \cdot \alpha^*(f).
\end{equation}

\medskip

\begin{prop} \label{Unitary}
The map $U$ extends 
to a unitary operator from 
$L^2( \alpha(I), ds)$ 
onto $L^2(I, \rho^{n-2}dt)$.
Moreover,
\begin{equation} \label{Conjugated}
U^{-1} \circ L \circ U~=~
   \frac{\partial^2}{\partial s^2}~ 
  - \frac{bd}{2} \cdot (\alpha^{-1})^*
      \left( \rho^{-2a-2} \left(
   \rho \cdot \rho_{tt}''
        +\frac{bd-2a-2}{2}(\rho'_t)^2 \right)
   \right).
\end{equation}
\end{prop}
\begin{proof}
The first claim is straightforward.
The second claim is a lengthy but
straightforward computation.
\end{proof}

The fibrewise 
$L^2(\{t\} \times M, dV_h)$-norm of a 
function $f \in C^0(I \times M)$ is 
the function on $I$ defined by
\begin{equation} 
   ||f||^2_M (t) = \int_{\{t\}\times M}
    f^2(t,m) dV_h(m).
\end{equation}

\medskip

\begin{prop} \label{BasicProp}
Let $\psi \in C^{2}(I \times M)$
satisfy $\Delta_g \psi =\lambda \psi$.
Then
\begin{equation} \label{LDiff}
  \frac{1}{2} L \left( ||\psi||^2_M \right)~
  = -\lambda \cdot ||\psi||^2_M~ +
     \int_{\{t\} \times M} 
     g(\nabla \psi, \nabla \psi)~ dV_h.
\end{equation}
\end{prop}

\medskip

\begin{proof}
Straightforward computation gives
\begin{equation*} 
  \frac{1}{2} L ( ||\psi||^2_M )~ 
  =~ \int_{M} L (\psi^2) dV_h~ 
  =~   \int_{M} \psi L (\psi) dV_h 
  + \rho^{-2a} \int_{M}
 (\partial_t \psi)^2~ dV_h. 
\end{equation*}
By hypothesis, we have 
$\int \psi \Delta \psi = \lambda \int \psi^2$
and hence by (\ref{Laplacian})
\begin{equation*} \label{Energy}
  -\int_{M} \psi L (\psi)~ dV_h
  + \rho^{-2b} \int_{M}
 \psi \cdot \Delta_h \psi~ dV_h~
  =~ \lambda \int_M \psi^2~ dV_h. 
\end{equation*}
Integrating by parts over $M$
gives
\begin{equation}
\int_{M}\psi \cdot \Delta_h \psi~ dV_h
  = \int_{M}
    h(\nabla_h \psi, \nabla_h \psi)~ dV_h.
\end{equation}
Also note that from (\ref{Gradient})
we have 
\begin{equation} \label{GDecomposed}
  g(\nabla \psi, \nabla \psi) =
  \rho^{-2a} (\partial_t \psi)^2
  + \rho^{-2b}  
  h(\nabla_h \psi, \nabla_h \psi).
\end{equation}
The claim follows.
\end{proof}

\medskip 

Let $\mu_1$ denote the smallest 
non-zero eigenvalue of $\Delta_h$.

\medskip

\begin{lem}[Convexity] \label{ConvexLemma}
For any Laplace eigenfunction $\psi$ on
$I \times M$ having constant mode
$\psi_0\equiv 0$,
\begin{equation}
  L ( ||\psi||^2_M ) 
  \geq 2 \cdot \left(\mu_1 \cdot 
  \rho^{-2b} -\lambda \right) \cdot
  ||\psi||^2_M.
\end{equation}
\end{lem}

\medskip

\begin{proof}
Since $\psi_0 \equiv 0$, 
the  function $m \rightarrow \psi(t,m)$
is orthogonal to the constants. 
Since the $0$-eigenspace consists of the
constants, by the minimax principle
\begin{equation} \label{Poincare}
\int_{\{t\} \times M}
  h(\nabla_h \psi, \nabla_h \psi)~ dV_h~
 \geq~ \mu_1 \int_{\{t\} \times M} \psi^2 
 dV_h.
\end{equation}
The claim then follows from 
Proposition \ref{BasicProp}
and (\ref{GDecomposed}).
\end{proof}

\medskip

%%%%%%%%%%%%%%%%%%%%%%%%%%%%%%%%%%%%%%%%%%%%%

%% Manifolds with cusps

\section{The spectral theory of 
              manifolds with cusps}

\label{SectionManifoldCusps}

\begin{ass}
In this section, we assume 
$a \leq -1$, $b>0$, and $c>0$.
\end{ass}

The manifold $I \setminus \{0\} \times M$
equipped with the metric $g_0$ of 
(\ref{LocalMetric}) is a disjoint union
of `cusps'. Here we use the $\epsilon=0$ case
of the analysis in the previous section 
to derive results about
manifolds with cusps.

Let $I^{+} \subset {\mathbb R}^+$ be an 
interval with lower endpoint equal to zero.
The manifold 
$I^{+} \times M$ equipped with the metric
\begin{equation} \label{Cusp}
 g~ =~ (ct)^{-2a}~ dt^2~ + (ct)^{2b}~ h.
\end{equation}
will be called a {\em cusp} of type
$(a,b,c)$. 
Note that the limiting manifold 
$(N^0 ,g_0)$ of an 
$(a,b)$ degenerating 
family is a manifold with ends 
isometric to cusps.

In the following, we let $\sigma_{ess}(T)$ 
denote the essential spectrum of an operator 
$T$.
\medskip

\begin{prop} \label{ContinuousSpectrum}
Let $(N,g)$ be a $d$-dimensional
Riemannian manifold with finitely 
many ends each of which is a cusp
of type $(a_j, b_j, c_j)$. 
If for some $j$, we have $a_j<-1$, then 
$\sigma_{ess}(\Delta_g)=[0, \infty[$.
Otherwise each $a_j=-1$, and
we have $\sigma_{ess}(\Delta_{g})= [ m , \infty[$
where $m$ is the infimum of
\begin{equation} \label{Bottom}
 \left(\frac{c_j \cdot 
          b_j \cdot d_j}{2} \right)^2.
\end{equation}
\end{prop}
\medskip 

\begin{proof}
By assumption $a \leq -1$, and hence 
$(N,g)$ is Riemannian complete \cite{Mlr}
\S 8.1. 
Therefore $\sigma_{ess}(\Delta)$ depends only 
on the geometry of the ends. (See, 
for example, \cite{DnlLi79} Proposition 2.1).
In particular, let $\Delta_j$ be the 
Friedrichs extension of $\Delta$ 
restricted to smooth functions supported 
in the $j^{th}$ end with respect to 
the $L^2$-norm induced by $dV_g$. 
Then
\begin{equation}
 \sigma_{ess}(\Delta)~ =~ 
 \cup_{j}~ \sigma_{ess}( \Delta_j). 
\end{equation}
Hence the claim reduces to the 
consideration of a single cusp.
From (\ref{Diagonal}) we have that 
\begin{equation}
  \sigma_{ess}(\Delta)~ =~ \sigma_{ess}(-L)
  \cup \sigma_{ess}(\Delta^{\perp}).
\end{equation}
Hence the claim is a consequence of 
Propositions \ref{DeltaPerp}
and \ref{LSpectrum} below.
\end{proof}

\begin{prop}  \label{LSpectrum}
If $a=-1$, then $\sigma_{ess}(-L)
    = [(2^{-1} \cdot c \cdot b\cdot d)^2, \infty[$.
If $a<-1$, then 
$\sigma_{ess}(-L)= [0, \infty[$.
\end{prop}
\begin{proof}
Let $U$ be as in Proposition \ref{Unitary}
where $\epsilon=0$. Since $U$ is unitary, 
we have $\sigma_{ess}(-L)
    = \sigma(-U^{-1} \circ L \circ U)$.
Moreover, since $\rho(0,t) = c t$,
$\rho''_{tt}=0$  and $(\rho'_t)^2=c^2$.
Thus, with $a=-1$  and $\epsilon=0$,
identity (\ref{Conjugated}) specializes to
\begin{equation} \label{ConjugateZero}
  U^{-1} \circ L \circ U~
  =~ \partial_s^2 - 
\left(\frac{ c \cdot b \cdot d}{2} \right)^2.
\end{equation}
The first claim follows.

If $a<-1$, then (\ref{Conjugated}) becomes
\begin{equation} 
  U^{-1} \circ L \circ U~
  =~ \partial_s^2 -  
  k \cdot (\alpha^{-1})^*((ct)^{-2a-2})(s).
\end{equation}
where $k= 4^{-1} \cdot bd \cdot (bd-2a-2)$.
From (\ref{Alpha}) we have
$\alpha(t) \sim (a+1)^{-1} 
           \cdot c^a \cdot t^{a+1}$
as $t \rightarrow 0$. Thus,
\begin{equation}
  (\alpha^{-1})^*((ct)^{-2a-2})(s) \sim  
     (c \cdot (a+1) \cdot s)^{-2}
\end{equation}
as $s$ tends to $-\infty$. Since $s^{-2}$
belongs to $L^2(]-\infty, -1], ds)$, the
second claim follows from standard results
on the essential spectrum of Schr\"odinger
operators. See, for example, Theorem XIII.15
\cite{RdSmn}.
\end{proof}

\medskip

The proof of the following is modeled
on the proof of Lemma 8.7 in \cite{LaxPhl}.

\medskip

\begin{prop} \label{DeltaPerp}
The operator $\Delta^{\perp}$ 
is compactly resolved. Hence 
$\sigma_{ess}(\Delta^{\perp})= \emptyset$.
\end{prop}
\begin{proof}
Without loss of generality, 
the homogeneity constant $c$ equals $1$.
Note that the operator $\Delta^{\perp}$ 
is a Friedrichs extension
associated to the form 
\begin{equation} \label{Q}
 Q(f)~ =~ \int_{I^+ \times M}
 \left(t^{-2a} \cdot (\partial_t f)^2 +  
 t^{-2b} \cdot  h(\nabla_h f, \nabla_h f) 
 \right)~ dV_g.
\end{equation}
In particular, to prove the claim 
it will suffice to show 
that the intersection of the $Q$-unit 
ball with $L^2_{\perp}(I^+ \times M, dV_g)$
is compact in 
$L^2_{\perp}(I^+ \times M, dV_g)$. 
(See, for example,
Theorem XIII.64 \cite{RdSmn}).

We claim that for any 
$\delta \in I^+=[0,t_0]$ and
any (smooth) function $\phi$ with 
$Q(\phi, \phi) \leq 1$ and $\phi_0 \equiv 0$,
we have 
\begin{equation} \label{OutInCusp}
 \int_{[0, \delta] \times M}
 \phi^2 dV_g~ \leq~ \frac{\delta^{2b}}{\mu_1}   
\end{equation}
where $\mu_1$ is the smallest nonzero eigenvalue
of $\Delta_h$.
Indeed, (\ref{Poincare}) holds and hence
\begin{eqnarray}
 \mu_1 \cdot
 \int_{[0, \delta] \times M} \phi^2~ dV_g~
 &\leq&
  \int_{[0, \delta] \times M} 
 h(\nabla_h~ \phi,
 \nabla_h~ \phi)~ dV_g  \\
 &=&  \delta^{2b}  \int_{[0, \delta] \times M} 
 t^{-2b} \cdot h(\nabla_h~ \phi,
 \nabla_h~ \phi)~ dV_g.  
\end{eqnarray}
The claim then follows 
from inspecting (\ref{Q}).

Now let $f_i$ be a sequence of functions
such that $Q(f_i, f_i) \leq 1$.
Since $b>0$, given $\epsilon>0$ small, 
there exists $\delta <t_0$ such that
$\delta^{2b}/\mu_1= \epsilon/2$. 
By Rellich's theorem, the restriction of
$f_i$ restricted to $[\delta,t_0] \times M$ 
is an $L^2$ Cauchy sequence. In particular,
there exists $N$ such that if $i,j>N$,
then 
\begin{equation} 
 \int_{[\delta,t_0] \times M}
 (f_i-f_j)^2~ dV_g~  \leq~ \frac{\epsilon}{2}.
\end{equation}
Applying (\ref{OutInCusp}) 
to $\phi=f_i-f_j$  then gives the claim.
\end{proof}

\medskip

The following Proposition gives the 
decay of `cusp forms'.

\medskip

\begin{prop} \label{CuspDecay}
Let $(N,g)$ be a Riemannian manifold with 
a cusp $(I \times M, g)$.
Let $\psi \in L^2(N, dV_g)$ 
be a Laplace eigenfunction
of $\Delta_{g}$ with
$\psi_0 \equiv 0$ on $I \times M$.
Then $||\psi||^2_M$ is convex
and for all $j>0$
\begin{equation} \label{TDecay}
  \lim_{t \rightarrow 0}
  t^{-j} \cdot ||\psi||^2_M(t) = 0.
\end{equation}
\end{prop}

\medskip

\begin{proof}
Since $b>0$, by Lemma \ref{ConvexLemma}
we have that for any $k>0$, there exists
a $t_*>0$ such that for all $t \in [0,t_*]$
\begin{equation} \label{ConvexAp}
  L(||\psi||^2_{M})~ \geq~ k^2 ||\psi||^2_{M}.
\end{equation}
For $U$ defined as in (\ref{UnitaryDefn}),
let $f= U^{-1}( ||\psi||^2_{M})$. 
Note that we have 
\begin{equation} \label{Mono}
 h \geq g~ \mbox{ if and only if }~ U(h) 
 \geq U(g).
\end{equation}
It follows that (\ref{ConvexAp}) holds iff
\begin{equation} \label{ConvexAp2}
  U^{-1} \circ L \circ U (f)~
  \geq~ k^2 f.
\end{equation}
Since  $\psi \in L^2(N, dV_g)$ and
$U$ is a unitary transformation, we have 
$f \in L^2(]-\infty, \alpha(t_0)], ds)$.
Moreover, since $b>0$ and $a \leq -1$, 
from (\ref{Conjugated}) we have
\begin{equation}
 U^{-1} \circ L \circ U f= \partial_s^2 f - \phi \cdot f
\end{equation}
where $\phi \geq 0$. Thus, it follows from (\ref{ConvexAp2})
that $f$ is convex and that $f(s)=O(e^{ks})$  as $s \rightarrow -\infty$. 

Note that if $h$ is a linear, then
since $b>0$, and $a \leq -1$, 
the function $U(h)$ is also convex. 
Thus, since $f$
is conxex, it follows from (\ref{Mono}) 
that $U(f)=||\psi||^2_M$ is convex.
Moreover, from the definition of $U$,
we find that $U(f)=||\psi||^2_M$ 
satisfies (\ref{TDecay}) with 
$j= k-\frac{bd}{2}$. 
\end{proof}

\medskip

\begin{lem} \label{Decay}
Let $(N,g)$ be a manifold with a
cusp $(I \times M,g)$.
Let $\psi \in L^2(N, dV_g)$ 
be a Laplace eigenfunction
of $\Delta_{g}$ with
$\psi_0 \equiv 0$ on $I \times M$.
Then $\partial_t \psi$ belongs
to $L^2(I \times M, dV_g)$.
\end{lem}
\begin{proof}
Since $b>0$, there exists
$t_*>0$ such that for $0<t<t_*$, we have 
$\lambda -\mu_1 t^{-2b} < 0$. 
It then follows 
from Proposition \ref{BasicProp}, 
(\ref{Poincare}), and (\ref{GDecomposed}), 
that for $0<t<t_*$
\begin{equation}
 (ct)^{-2a}  \int_M (\partial_t \psi)^2~
 \leq~ \frac{1}{2} \cdot 
 L(||\psi||^2_M).
\end{equation}
Multiply both sides by $t^{2a}$
and integrate over $I_s=[s,t_*]$ to obtain
\begin{equation}
 \int_{I_s \times M} 
 (\partial_t \psi)^2~ dV_g
 \leq~ \frac{c^{2a}}{2} \cdot   
 \int_{s}^{t_*} t^{2a} \cdot L(||\psi||^2_M)~
 t^{a+bd}~dt.
\end{equation}
Thus, it suffices to show that the 
integral on the right hand
side remains bounded as $t$ tends to zero.
To verify this, we integrate by parts 
and obtain
\begin{equation} \label{Party}
\left. 
 \int_{s}^{t_*} t^{k} \cdot L(||\psi||^2_M)~
=~ 
 \left( t^{\alpha} \cdot \partial_t ||\psi||^2_M 
+ C_1 \cdot t^{\beta} ||\psi||^2_M \right) \right|_{s}^{t_*}~ +~
  C_2 \cdot \int_{s}^{t_*}
 t^{\gamma} \cdot  ||\psi||^2_M.
\end{equation}
where $C_1$, $C_2$, $\alpha$, $\beta$, and $\gamma$ are constants that depend on $k$.
By Proposition \ref{CuspDecay}, $||\psi||^2_M$ is convex and satisfies
(\ref{TDecay}). It follows that both $t^{\beta} ||\psi||^2_M$ 
and $t^{\gamma} \partial_t ||\psi||^2_M$ are bounded for sufficently small $t$.
The claim follows.
\end{proof}

%%%%%%%%%%%%%%%%%%%%%%%%%%%%%%%%%%%%%%%%%%%%

%%  CONVERGENCE OF EIGENFUNCTIONS SECTION

\section{Convergence of eigenfunctions}

\label{SectionEigenfunctionConvergence}

\medskip

\begin{thm}  \label{EigenfunctionConvergence}
Let $(N, g_{\epsilon})$ be a 
$(a,b)$-degenerating family with $a \leq -1$
and $b>0$. For $\epsilon_j \rightarrow 0$, 
let $\psi_{j}$ be a sequence of
eigenfunctions of $\Delta_{\epsilon_j}$ 
on $L^2(N, dV_g)$ with eigenvalues 
$\lambda_j \leq \Lambda$. Then there 
is a subsequence $\psi_k$, a sequence 
$a_k \in {\mathbb R}$, and a nontrivial 
eigenfunction $\psi_*$ of the $C^{\infty}$ 
Laplacian on  $N^0$, such that for every 
compact $A \subset N^0$,
the sequence $a_k \psi_k$ converges 
to $\psi_*$ in $H^1(A, dV_{g_0})$.
Moreover, if  $(\psi_k)_0 \equiv 0$ 
for each $m$, then $a_k \psi_k$ converges 
to $\psi_*$ in 
$L^2_{\perp}(N^0, dV_{g_0})$.
\end{thm}
\begin{proof}
For $f \in C^0(N)$, define 
\begin{equation} \label{Normal}
\overline{f}(x) 
  = \left\{ \begin{array}{ll} f(x)
  & x \in
  N \setminus (\frac{1}{3} I \times M)  \\
  f(x)-f_0(x)
  & x \in \frac{1}{3} I \times M \\
  \end{array}  \right\}
\end{equation}
where $f_0$ is the constant mode of $f$.
Now define  
\begin{equation} \label{Normal2}
 \phi_j 
 =  \frac{ \psi_j}{||\overline{\psi}_j||}
\end{equation}
where $||\cdot||$ is the ${L^2(N, dV_g)}$
norm. We will show that this renormalized
sequence satisfies the claim.

Let $\chi \in C^{\infty}_0(N)$ be nonnegative,
have support in 
$N \setminus (\frac{1}{3} I \times M)$, and
satisfy $\chi \equiv 1$ 
on $N \setminus (\frac{2}{3} \times M)$.
Using the Schwarz inequality (for both
$TM$ and $L^2$) and the fact that 
$Q(\phi_j)= \lambda_j \cdot ||\phi_j||$, 
we find that 
\begin{equation}
  Q(\chi \cdot \phi_j)~ \leq~
  \left( \left(\int_N \chi^2 |\nabla \phi_j|^2 
  dV_g \right)^{\frac{1}{2}}
 +  \left(\int_N \phi^2 |\nabla \chi|^2
  dV_g \right)^{\frac{1}{2}}  \right)^2
\end{equation}
Using integration by parts and the fact
that $\Delta \phi_j =\lambda_j \cdot \phi_j$, 
one finds a constant $C$ depending only on $\chi$
such that $\int_N \chi^2 |\nabla \phi_j|^2$
is less than $C \cdot (\Lambda +C)$ times 
the integral of $\phi_j^2$ over 
$N \setminus (\frac{1}{3} I \times M)$.  
By (\ref{Normal}),
this latter integral is less than 1, and hence
$Q(\chi \cdot \phi_j)$ is a bounded sequence
in the space of $H^1(N, dV_g)$ functions 
that have support outside 
of $\frac{1}{3} I \times M$. 
By Rellich's theorem we may pass 
to a subsequence such that 
$\chi \cdot \phi_j$  converges in this 
space to a function $\chi \cdot \phi_*$. 
Hence the restriction of $\phi_j$ to 
$N \setminus (\frac{2}{3} I \times M)$
converges in $H^1$.

It follows that the restriction of the  
constant mode, $(\phi_j)_0$, to
$I \setminus (\frac{2}{3} I )$ converges
to $(\phi_*)_0$ 
in $H^1$ norm. Thus, by Sobolev's
embedding theorem (in dimension 1)  
$(\phi_j)_0$ converges to
$(\phi_*)_0$ in $C^0$, and hence 
the boundary conditions 
of (\ref{ConstantODE}) on $I \setminus K$ 
converge in $C^0$ as $j \rightarrow \infty$
for each nontrivial interval 
$0 \in K \subset I$. Since the 
metrics converge, the coefficients of 
(\ref{ConstantODE}) converge  
in $C^0$. It follows that $(\phi_*)_0$ 
extends uniquely to $I\setminus \{0\}$, 
and that, for each $K \ni 0$, the sequence
$(\phi_j)_0$ converges to $(\phi_*)_0$
in $C^0(K \setminus I)$.

Since $\phi_j$ is obtained from 
$\overline{\phi}_j$ by adding on the 
constant mode $(\phi_j)_0$, we find that
for every interval $K \subset I$,
the restriction of $\phi_j$ to 
$N \setminus (K \times M)$ is
uniformly bounded in $L^2$-norm.
Hence, for each $K \subset I$, the 
argument above can be applied to give
a further subsequence $\phi_j$  
whose restriction converges in 
$H^1(N \setminus (K \times M), dV_g)$.
Diagonalization yields a further 
subsequence
that converges to a function $\phi_*$ 
for every $K$.

Since $\lambda_j$ is bounded, we may take
a further subsequence such that 
$\lambda_j$ converges to some 
$\lambda_* \geq 0$. For each test
function $T$ supported away from 
$\{0\} \times M$, we have that  
$Q_{\epsilon_j}(\phi_j, T) 
\rightarrow Q_0(\phi_*, T)$. It follows
that $\phi_*$ is a weak---and hence 
by elliptic regularity, a 
strong---solution to 
$\Delta_0 \phi_* = \lambda_* \cdot \phi_*$.

It remains to show that $\phi_*$ does not
vanish identically.  Since $b >0$ and $a \leq -1$, 
there exists an interval $J \subset \frac{1}{3} I$ 
symmetric about $0$  such that for all $t \in J$
and $\epsilon$ sufficiently small
\begin{equation} \label{Hi}
\rho(\epsilon,t)^{2b}~ \leq~ 
     \frac{2 \mu_1}{\Lambda}
\end{equation}
and 
\begin{equation} \label{Led}
 L(\rho^{2b})(\epsilon,t)~ \leq~ \frac{\mu_1}{2}.
\end{equation}
Indeed, the operator $L$ 
adds $-2a-2\geq 0$ to the degree of
homogeneity. Let $\chi \in C_{0}^{\infty}(J)$
with $\chi \equiv 1$ on $\frac{1}{2} J$. 
From (\ref{Hi}), Lemma \ref{ConvexLemma}, and 
the self-adjointness of $L$, we obtain:
\begin{eqnarray} \nonumber
\mu_1 \cdot \int_{J \times M} 
 \chi \cdot \overline{\phi}_j^2~ dV_g~ 
 &=& 2 \int_{J}  \chi \cdot
 \rho^{2b} \cdot 
 (\mu_1 \rho^{-2b} -\lambda_j) ||\overline{\phi}_j||^2_M~
 \rho^{a+bd}~ dt  \\ \nonumber
 &\leq& \int_{J}  \chi \cdot
 \rho^{2b} \cdot L(||\overline{\phi}_j||^2_M)~
 \rho^{a+bd}~ dt  \\
 &=& \int_{J \times M} L( \chi \cdot
   \rho^{2b}) \cdot \overline{\phi}_j^2~
  dV_g.  \label{AfterPart}
\end{eqnarray}
Note that $L( \chi \cdot \rho^{2b})$ equals
$\chi \cdot L(\rho)$ plus a smooth
function supported on 
$J \setminus \frac{1}{2} J$ bounded by 
$C$. Thus, it follows from 
(\ref{Led}) and (\ref{AfterPart}) that
\begin{equation}
 \frac{\mu_1}{2} 
\int_{\frac{1}{2}J \times M} 
\overline{\phi}_j^2~ dV_g~
 \leq~  C \cdot 
 \int_{(J\setminus \frac{1}{2} J) \times M}
 \overline{\phi}_j^2~ dV_g.
\end{equation}
Thus, from the normalization of 
the $L^2$-norm in (\ref{Normal2}), 
we find that 
\begin{equation}
1~ \leq ~\frac{2 C}{\mu_1}
 \int_{(J \setminus \frac{1}{2} J) 
\times M}
  \overline{\phi}^2~ dV_g~ +~
   \int_{N \setminus 
   (\frac{1}{2}J \times M)}
 \overline{\phi}^2~
 dV_g.
\end{equation}
Therefore, since $\phi_j$ restricted to 
$N \setminus ( \frac{1}{2} J \times M)$
converges in $L^2$ to $\phi_*$, 
the function $\phi_*$ is nontrivial.
\end{proof}

%%%%%%%%%%%%%%%%%%%%%%%%%%%%%%%%%%%%%

% Delta^{perp} spectral continuity

\section{Eigenvalue  continuity for 
                    $\Delta^{\perp}$}

\label{SectionCutOffLaplacian}

In this section, $g_{\epsilon}$ 
is a family of $(a,b)$-degenerating 
metrics on $I \times M$ that depends
continuously on $\epsilon$.  
By Proposition \ref{DeltaPerp},
for each $\epsilon$---including 
$\epsilon=0$---the operator 
$\Delta^{\perp}_{\epsilon}$
has a countable collection of eigenvalues
(including multiplicities)
\begin{equation}
0 < \lambda_1(\epsilon) 
   \leq \lambda_2(\epsilon) 
  \leq \lambda_3(\epsilon) \leq \cdots.
\end{equation}

\medskip

\begin{thm} \label{SpectralContinuity} 
Let $a \leq -1$ and $b>0$. 
Then for each $i$, the function
$\epsilon \rightarrow \lambda_i(\epsilon)$ 
is continuous.
\end{thm}

\medskip

\begin{proof}
The operator $\Delta^{\perp}_{\epsilon}$ 
is defined via the (sesquilinear) form
\begin{equation} \label{QuadraticForm}
  Q_{\epsilon}(f)= \int \rho^{-2a} 
 \cdot 
 (\partial_t f)^2 \cdot \rho^{a+bd}~ dV_h 
  + \int \rho^{-2b} \cdot
  h(\nabla_h f, \nabla_h f) \cdot
  \rho^{a+bd}~dV_h
\end{equation}
restricted to 
$C^{\infty}_0(I \times M) 
     \cap L^2_{\perp}(I 
         \times M, dV_{g(\epsilon)})$.
For $\epsilon \neq 0$, the domain of
$\Delta^{\perp}_{\epsilon}$ is 
$H^1(I \times M, dV_g) \cap L^2_{\perp}(I 
         \times M, dV_{g(\epsilon)})$.  
For any power $c$, the function $\rho^{c}$ 
is uniformly continuous
in $\epsilon$ over closed intervals
that do not contain $0$. It follows from 
\S VI.3.2 \cite{Kat}
that $\Delta_{\epsilon}^{\perp}$ is continuous
in the `generalized sense' 
for $\epsilon \neq 0$.
Thus, by \S IV.3.5 \cite{Kat},
any finite system of eigenvalues varies
continuously for $\epsilon \neq 0$. 
It follows that each 
$\lambda_i$ is continuous for 
$\epsilon \neq 0$.

We are left with showing the continuity
at $\epsilon=0$. 
Let $V_{k-1}(\epsilon)$ be
a span of eigenfunctions associated to 
the first $k$ eigenvalues:
$\lambda_1(\epsilon), 
\ldots, \lambda_{k}(\epsilon)$. 
Then by the minimax principle, 
$\lambda_{k+1}(\epsilon)$ 
is the minimum value 
of the functional $F_{\epsilon}(\phi)= Q_{\epsilon}(\phi)/ ||\phi||^2_{\epsilon}$ 
over $V_{k-1}(\epsilon)^{\perp}$ where
$||\cdot||_{\epsilon}$ 
denotes the 
$L^2(I \times M, dV_{g_{\epsilon}})$-norm.

Let $\phi_1$ be an eigenfunction of 
$\Delta^{\perp}_0$ with eigenvalue 
$\lambda_1(0)$.
By Proposition \ref{Decay}, 
the function $\phi_1$ belongs to 
the domain of 
$F_{\epsilon}$ for small $\epsilon$.  
From the continuity of $\rho$, we find 
that 
\begin{equation*}
 \lim_{\epsilon \rightarrow 0} 
 F_{\epsilon}(\phi_1) =
 F_{0}(\phi_1) = \lambda_1(0).
\end{equation*}
Hence, by the minimax principle, 
$\limsup_{\epsilon \rightarrow 0} 
   \lambda_1(\epsilon) \leq \lambda_1(0)$.
For each $\epsilon$, let $\phi_1(\epsilon)$ 
be an eigenfunction of 
$\Delta^{\perp}_{\epsilon}$ with 
eigenvalue $\lambda_1(\epsilon)$.
By applying Theorem 
\ref{EigenfunctionConvergence} to a 
subsequence whose eigenvalues limit 
to $\liminf_{\epsilon \rightarrow 0}
            \lambda_1(\epsilon)$, 
we obtain an eigenfunction $\phi$ of 
$\Delta^{\perp}_{0}$ with eigenvalue 
$\lambda_*=\liminf_{\epsilon \rightarrow 0} 
     \lambda_1(\epsilon)$ 
less than or equal to $\lambda_1$. 
But since $\lambda_1$ is the smallest 
eigenvalue, $\lambda_*=\lambda_1$.  
It follows that
$\lim_{\epsilon \rightarrow 0} 
      \lambda_1(\epsilon)= \lambda_1(0)$.

The continuity at $\epsilon=0$ 
of $\lambda_k$ for general $k$ follows from 
a straightforward inductive argument
involving $V_{k}(\epsilon)$. (Note that no 
claim is made about the continuity of 
the family $V_{k}(\epsilon)$.)
\end{proof}

\begin{remk}
With minor modifications, the
proof of Theorem \ref{SpectralContinuity} 
gives the continuity of the eigenvalues
of the {\em cut-off (or pseudo-)Laplacian} 
defined as in \cite{Ji93}.
\end{remk}

%%%%%%%%%%%%%%%%%%%%%%%%%%%%%%%%%%%%%%%%%%%%%%%%%

\section{Counting 
           relatively small eigenvalues}

\medskip

\label{SectionCount}

\medskip

In the following $O(1)$ will denote a bounded function of $\epsilon$. 
\medskip

\begin{thm} \label{CountTheorem}
Let $a=-1$ and $b>0$.
For $\Lambda>0$, let 
${\mathcal N}_{\Lambda}(\epsilon)$ 
be the number of eigenvalues of 
$\Delta_{g_{\epsilon}}$ that lie 
in $[0, \Lambda]$. Then 
\begin{equation}  \label{MainEigenvalueCount}
  {\mathcal N}_{\Lambda}(\epsilon) =
   \left( c_+
 \sqrt{\Lambda-
  \left(\frac{c_{+} \cdot b \cdot d}{2}
  \right)^2 } 
  +   c_{-} \sqrt{\Lambda- 
 \left( \frac{c_{-}\cdot b \cdot d}{2}
  \right)^2 }  \right) \cdot
   \frac{\log( \epsilon^{-1})}{\pi} 
   + O_{\Lambda}(1).
\end{equation}
Here the value of the square root is taken to be zero if the argument
is negative.
\end{thm}
\begin{proof}
By the Dirichlet monotonicity, we have that 
${\mathcal N}_{\Lambda}(\epsilon)$ is
 bounded below by the number, 
${\mathcal N}_{\Lambda}^D(\epsilon)$, 
of eigenvalues of the Dirichlet problem 
on $I \times M$.  By Neumann monotonicity,  
${\mathcal N}_{\Lambda}(\epsilon)$ is 
bounded above by the number, 
${\mathcal N}_{\Lambda}^N(\epsilon)$, 
of eigenvalues of the Neumann problem on 
$I \times M$ plus the number of eigenvalues 
of the Neumann problem on 
$N \setminus (I \times M)$.
The latter is $O_{\Lambda}(1)$ since 
$g(\epsilon)$ converges uniformly
on $N \setminus (I \times M)$. In sum, 
we have
\begin{equation}
   {\mathcal N}_{\Lambda}^{D}(\epsilon) \leq 
       {\mathcal N}_{\Lambda}(\epsilon) 
     \leq {\mathcal N}_{\Lambda}^{N}
    (\epsilon) + O_{\Lambda}(1).
\end{equation}

By (\ref{Diagonal}) the Dirichlet
(resp. Neumann)  spectrum decomposes 
into the Dirichlet (resp. Neumann) 
eigenvalues of $L$ acting on 
$L^2(I, \rho^{n-2}dt)$
and those of $\Delta^{\perp}_{\epsilon}$ 
acting on 
$L^2_{\perp}(I \times M, dV_{g_{\epsilon}})$.
By Theorem \ref{SpectralContinuity}, 
the number of Dirichlet (resp. Neumann) 
eigenvalues of $\Delta^{\perp}_{\epsilon}$
is $O_{\Lambda}(1)$.  Hence, the 
claim reduces to the Lemma \ref{CountODE}
below. 
\end{proof}

\medskip

To emphasize its dependence on $\epsilon$,
we let $L_{\epsilon}$ denote the operator
defined in (\ref{L}).

\medskip

\begin{lem} \label{CountODE}
For each $\epsilon>0$, let 
${\mathcal N}^{\partial}_{\Lambda}(\epsilon)$
be the number of solutions 
$\lambda \in [0, \Lambda]$  
to the Dirichlet (resp. Neumann) 
eigenvalue problem 
\begin{equation}
 \label{BasicEigenvalueProblem}
 -L_{\epsilon} v= \lambda v
\end{equation}
on the interval $I$.  Then
${\mathcal N}^{\partial}_{\Lambda}(\epsilon)$
satisfies 
(\ref{MainEigenvalueCount}).
\end{lem}
\begin{proof}
First note that since $a=-1$ and  
$\rho$ is homogeneous of degree 1,
 the eigenvalue problem 
(\ref{BasicEigenvalueProblem}) 
on $I$ is equivalent to the problem 
\begin{equation} \label{AfterHomogeneity}
L_1 v= -\lambda v
\end{equation}
on the dilated interval $\epsilon^{-1} I$. 
By conjuagting both sides of 
(\ref{AfterHomogeneity}) by 
$U=U_1$ and applying Proposition 
\ref{Unitary} , we obtain
the equivalent eigenvalue problem
\begin{equation} \label{Modified}
\partial_t^2 u -\frac{bd}{2} 
(\alpha_1^{-1})^* \left( 
\rho \cdot \rho_{tt}''
 +\frac{bd}{2}(\rho'_t)^2 \right) u~ 
=~ -\lambda u.
\end{equation}
on the interval $\alpha_1(\epsilon^{-1} I)$.

Let $a_{-}(\epsilon) <  
    a_{+}(\epsilon)$ be the endpoints
of the interval $\alpha_1(\epsilon^{-1} I)$.
Let ${\mathcal N}^{+}_{\Lambda}(\epsilon)$
(resp. ${\mathcal N}^{-}_{\Lambda}(\epsilon)$)
be the number of
eigenvalues $\lambda$ of (\ref{Modified}) 
on the interval $[0,a_+]$ (resp. $[a_-,0]$).
Then by Dirichlet-Neumann bracketing 
(see, for example, \cite{CrnHlb} p 408-409), 
we have 
\begin{equation} \label{Split}
 {\mathcal N}^{\partial}_{\Lambda}(\epsilon)
  = {\mathcal N}^{+}_{\Lambda}(\epsilon)+
  {\mathcal N}^{-}_{\Lambda}(\epsilon)
  +O_{\Lambda}(1).
\end{equation}

We claim that it suffices to show that
\begin{equation} \label{AAsymp}
 a_{\pm}(\epsilon)
   \sim \pm c_{\pm} \ln(\epsilon^{-1}),
\end{equation}
as $\epsilon$ tends to zero, and that
\begin{equation} \label{Equivalence}
  \frac{bd}{2}(\alpha_1^{-1})^* \left( 
  \rho \cdot \rho_{tt}''
  +\frac{bd}{2}(\rho'_t)^2 \right)~
  =~ \left(
 \frac{c_{\pm} \cdot b \cdot d}{2}
  \right)^2 
    \cdot u  +r(s)
\end{equation}
where $r(s)$ is $O(|s|^{-2})$ for $|s|$ large. 
Indeed, then (\ref{Modified}) would be 
equivalent to
\begin{equation} \label{NormalODE}
 \partial_s^2 u~ =~ 
 \left(- \lambda  
 +\left( \frac{c_{\pm} \cdot b \cdot d}{2}
  \right)^2 
\right)u  +r(s).
\end{equation}
By Proposition \ref{ClassicalCount} 
below, we would then have
\begin{equation}
{\mathcal N}^{\pm}_{\Lambda}(\epsilon)
  =  a_{\pm}(\epsilon) \sqrt{
 \Lambda - \left(
 \frac{c_{\pm} \cdot b \cdot d}{2}
  \right)^2},
\end{equation}
and the claim would follow from (\ref{Split}).

To verify (\ref{AAsymp}) and 
(\ref{Equivalence}), we use
homogeneity. We have 
$\rho^{-1}(1,t) \sim c_{\pm} t^{-1}$,
and hence $\alpha_1(t) \sim c_{\pm} \ln(t)$ 
as $t \rightarrow \pm \infty$.
Estimate (\ref{AAsymp}) follows. 
Moreover, 
\begin{equation} \label{AlphaInverse}
  \alpha_1^{-1}(t) 
          \sim \exp(c_{\pm} t).
\end{equation}

To prove (\ref{Equivalence}), we use
the following fact: If $\sigma(\epsilon,t)$ 
is smooth, homogeneous of degree $k$,
and $\sigma(0,t) \equiv 0$, then 
$\sigma(1,t) = O(|t|^{-{k+1}})$ as 
$|t|$ tends to infinity. (Indeed,
we have $h(1,t)= t^{-k} h(t^{-1}, 1)$
and Taylor expand in the first coordinate.)
By applying this fact to 
$\rho''_{tt}$ and $(\rho'_t)^2-c_{\pm}^2$, we find that
\begin{equation}
\rho \cdot \rho_{tt}''(1,t) = O(|t|^{-1})
\end{equation}
and 
\begin{equation}
 (\rho'_t)^2  = c_{\pm}^2 + O(|t|^{-1}).
\end{equation}
Equation  (\ref{Equivalence}) then follows
from (\ref{AlphaInverse}).
\end{proof}

\medskip

The following appears in \cite{ChvDdz94}. 
To make our exposition self-contained,
we include it here with an alternate proof.

\medskip

\begin{prop}   \label{ClassicalCount}
Let $r \in C^{0}({\mathbb R})$ be 
integrable with respect 
to Lebesgue measure. Let 
${\mathcal N}(a,b)$ be the number 
of solutions $\mu \in [0, M]$ to 
the Dirichlet (resp. Neumann) eigenvalue 
problem
\begin{equation} \label{IdealProblem}
\partial_s^2 u = -\mu u  +r u
\end{equation}
 on $[0, a]$. Then
\begin{equation}
  {\mathcal N}_M(a) = \frac{a \sqrt{M}}{\pi} + O_{M}(1)
\end{equation}
where $O_{M}(1)$ is a bounded function of $a \in {\mathbb R}$. 
\end{prop}
\begin{proof}
Let $0=\mu_0(a)< \mu_1(a) < \mu_2(a)< \cdots$  denote the (necessarily
simple) Dirichlet (resp. Neumann) eigenvalues of (\ref{IdealProblem})
for the interval $[0,a]$.  By the standard theory, for each $k$, 
$\mu_k(a)$ is a decreasing function of $a>0$ 
with $\lim_{a \rightarrow 0} \mu_k(a) = + \infty$. 
It follows that 
\begin{equation}
   N_M(a)= \Card \{ k \in {\mathbb N}: \mu_k(b) = M \mbox{ and } b \leq a\}.
\end{equation}
See Figure \ref{EigenbranchCardinality}.

\begin{figure}
  \centering
  \psfrag{mu}{$\mu$}
  \psfrag{M}{$M$}
  \psfrag{a}{$a$}
  \psfrag{s}{$s$}
    \includegraphics[angle=0]{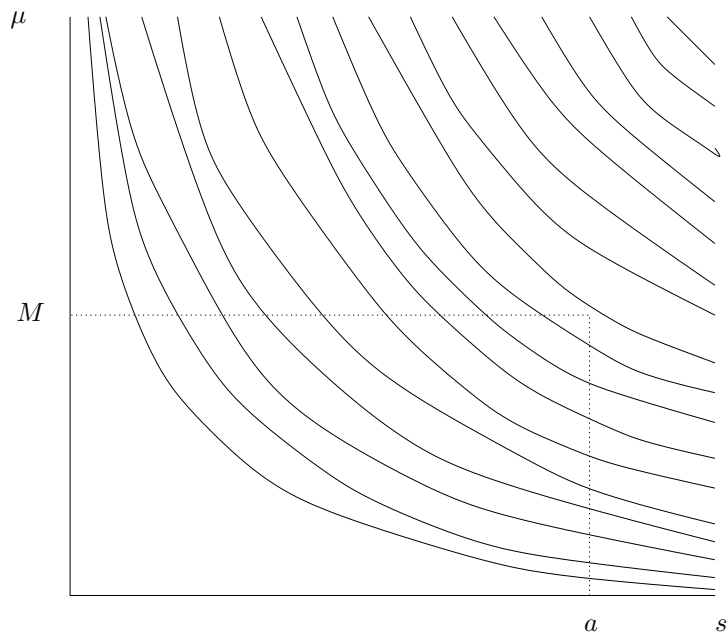}
  \caption{Sturm-Liouville eigenvalues on $[0,a]$ as functions of $a$.}
  \label{EigenbranchCardinality}
\end{figure}

Let $v \in C^2({\mathbb R})$ be a nonzero solution to (\ref{IdealProblem})
with $\mu=M$ and $v(0)=0$ (resp. $v'(0)=0$). (The function $v$ need not 
satisfy the Dirichlet (Neumann) boundary condition at $a$.). 
Since $r$ is integrable, there exist 
(see, for example, \cite{CrnHlb} p. 332-333) bounded functions
$\alpha, \delta \in C^1({\mathbb R})$, $\delta'(s)>0$, such that 
\begin{equation} \label{AsympBound} 
 v(s) = \alpha(s) \cdot \sin(\sqrt{M} s + \delta(s)).
\end{equation}
Note that $\mu_k(b)=M$ if and only if $v(b)=0$
(resp. $v'(b)=0$). Thus by (\ref{AsympBound}),
a Dirichlet eigenvalue $\mu_k(b)=M$ if and only if 
$\pi^{-1}(\sqrt{M} b + \delta(b)) \in {\mathbb Z}$
(resp. $\pi^{-1}(\sqrt{M} b + \delta(b)) \in {\mathbb Z} + \frac{\pi}{2}$).

In sum, for the Dirichlet problem,
\begin{equation}
  N_M(a) = \Card \{ 0\leq b \leq a: \pi^{-1}
            (\sqrt{M} b + \delta(b)) \in {\mathbb Z} \}.
\end{equation}
Thus, letting $C=\sup|\delta(s)|$ we have 
\begin{equation}
\pi^{-1} (\sqrt{M} b -C) -1  \leq 
   N_M(a)  \leq \pi^{-1}(\sqrt{M} b + C),
\end{equation}
and the claim follows for the Dirichlet case. 
To prove the Neumann case, one can use a similar analysis
involving differentiating (\ref{AsympBound}).  Or one can 
derive the Neumann case from the Dirichlet case via
Neumann-Dirichlet bracketing for Sturm-Liouville problems: 
the $k+2$nd Neumann eigenvalue is at least as great as the $k$th 
Dirichlet eigenvalue \cite{Wnb}.
\end{proof}

\end{document}